\documentclass[12pt]{article}

\normalbaselines
\catcode `\@=11
\@addtoreset{equation}{section}

\def\section{\@startsection {section}{1}{\z@}{-3.5ex plus -1ex minus
     -.2ex}{2.3ex plus .2ex}{\normalsize\bf}}
\def\subsection{\@startsection{subsection}{2}{\z@}{-3.25ex plus -1ex minus
 -.2ex}{1.5ex plus .2ex}{\normalsize\bf}}

\def\thebibliography#1{\section*{References\markboth
  {REFERENCES}{REFERENCES}}\list
  {[\arabic{enumi}]}{\settowidth\labelwidth{[#1]}\leftmargin\labelwidth
  \advance\leftmargin\labelsep
  \usecounter{enumi}}
  \def\newblock{\hskip .11em plus .33em minus -.07em}
  \sloppy
  \sfcode`\.=1000\relax}
 
\catcode `\@=12

\newtheorem{thm}{Theorem}
\newtheorem{rem}{Remark}
\newtheorem{cor}{Corollary}
\font\ac=eufm10 scaled\magstep1
\font\ab=msbm10 scaled\magstep1
\font\aba=msbm7 
\font\ai=eusb10 scaled \magstep 1
\newcommand{\gr}{\mbox{\ab R}}
\newcommand{\gcp}{\mbox{\ab CP}}
\newcommand{\gc}{\mbox{\ab C}}
\newcommand{\gcm}{\mbox{\aba C}}
\newcommand{\Gras}{\mbox{$G_n({\gc}^{m+n})$}}
\newcommand{\men}{\mbox{$\widetilde {\bf M}$}} 
\newcommand{\got}[1]{{\mbox{\ac{#1}}}}
\newcommand{\mb}[1]{{\mbox{\boldmath{$#1$}}}}
\newcommand{\gh}{\mbox{\ai H}}
\newcommand{\gl}{\mbox{\ai L}}
\newcommand{\gk}{\mbox{\ai K}}
\newcommand{\gph}{\mbox{{\bf P}(\gh )}}
\newcommand{\df}{\mbox{$:=$}}

\newcommand{\Ka}{K\"ahler}
\newcommand{\sfera}{\mbox{${\cal S}(\gh )$}}
\begin{document}
\vspace*{2.5cm}
\noindent
{ \bf  A REMARK ON BEREZIN'S QUANTIZATION AND CUT LOCUS}\\[1.3cm]
\noindent
\hspace*{1in}
\begin{minipage}{13cm}
Stefan Berceanu  \\[0.3cm]
Institute for Physics and Nuclear Engineering\\
Department of Theoretical Physics\\
PO BOX MG-6, Bucharest-Magurele, Romania\\
E-mail: Berceanu@theor1.ifa.ro; Berceanu@Roifa.Ifa.Ro\\


\end{minipage}

\begin{abstract}
\noindent
The consequences for Berezin's quantization on symmetric spaces of the identity
of the  set of  coherent vectors orthogonal to a fixed one with the cut locus
are stated precisely.  It is shown that functions expressing the coherent
states, the covariant symbols of operators,  the diastasis function, the
characteristic and two-point functions are defined  when one variable does
not belong to the cut locus of the other one.

\end{abstract}

\section{\hspace{-4mm}.\hspace{2mm} INTRODUCTION}

Recently I have emphasized \cite{ber2} the deep relationship between geodesics
and coherent states \cite{klauder}. For homogeneous manifolds in which the
exponential from the Lie algebra to the Lie group equals the geodesic
exponential, and in particular for symmetric spaces, it was proved that the cut
locus (${\bf CL}$) of a fixed point in the manifold is equal to the
so called
{\it polar divisor} (denoted $\Sigma$), i.e. the set of coherent vectors
orthogonal to the coherent vector corresponding to the fixed point.
The coherent states are  a powerful tool in global
differential geometry and algebraic geometry \cite{ber1}.

In this talk I comment the physical relevance of the result ${\bf CL}_0 =
\Sigma_0$. The coherent states offer a straightforward recipe
\cite{raw,cval,od} for geometric quantization \cite{kost}. Mainly, I point out
the consequences of the results obtained in Ref. \cite{ber2} for Berezin's
quantization, especially for symmetric spaces (see firstly the second reference
in \cite{berezin}). Also I clear up the significance
of  the identity ${\bf CL}_0 = \Sigma_0$    for the works of
Onofri \cite{onofri}, Rawnsley \cite{raw}, Moreno \cite{moreno} and Cahen, Gutt
and Rawnsley \cite{cgr}.

The main message of the present contribution is that on
symmetric spaces {\it the cut locus is present everywhere one speaks about
coherent states\/}. Indeed, the functions expressing the
coherent states, the covariant symbols of operators,
the diastasis function, the characteristic and two-point functions
are defined at a fixed $x$ for $y$ not in ${\bf CL}_x$.

The paper is laid out as follows: some mathematical results
are collected
in \S 2.  The notation in \S 3 for the coherent states differs from
that used earlier in Refs.
\cite{sbcag,ber7}. The main results of Ref. \cite{ber2} on the relationship
between coherent states and cut locus are briefly recalled in
\S 4, while the last Section is devoted to the announced comments.
I have also included in \S 5 a remark on the significance of the polar
divisor as divisor  in the meaning  of
algebraic geometry \cite{gh}, not exhibited in my talk, which
was generated by a\nopagebreak~ discussion during the Workshop. An
illustration of\nopagebreak~  all  presented results on the complex
Grassmann manifold \Gras~
is available in Ref. \cite{ber4}.
\section{\hspace{-4mm}.\hspace{2mm} MATHEMATICAL PRELIMINARIES}
1). Let $\chi$ be a  representation of the group $K$ on the Hilbert space
\gk~ and let us consider the principal bundle
\begin{equation}
\label{pbundle}
K\stackrel{i}{\rightarrow} G\stackrel{\lambda}{\rightarrow}\men,
\end{equation}
where \men~ is diffeomorphic with $G/K$, $i$ is the inclusion and
$\lambda $ is the natural projection $\lambda (g)=gK$. We
recall \cite{bott} the definition of the
$G$-{\it homogeneous vector bundle\/} ${\bf M}_{\chi}$ associated by the
character $\chi$
to the principal $K$-bundle (\ref{pbundle}), ${\bf M}_{\chi}\df
\men \times_{\chi}\gk$, or simply ${\bf M}\df\men \times_K\gk$. The
total space of ${\bf M}_{\chi}$ consists of all  equivalence
classes $[g,l]$ of elements $(g,l)$ under the
equivalence relation $(gp,l)\sim (g,\chi (p)l),~g\in G,~p\in K,~l \in
\gk$. If $U\subset\men$~ is open, let the notation
\begin{equation}
\label{local}
(G)^{U}=\{g\in G|\pi (g)\psi_0\in U\},
\end{equation}
where $\pi$ is a representation of $G$ whose restriction to $K$ is
$\chi$ and $\psi_0\in \gk$ corresponds to the base point $o\in \men$.
Then the continuous sections of
${\bf M}_{\chi}$ over $U$ are precisely the continuous maps
$\sigma :U\rightarrow G\times_{\chi}\gk$ of the form
\begin{equation}
\label{sectiuni}
\sigma(\pi (g)\psi_0)=[g,e_{\sigma}(g)],~e_{\sigma}:(G)^U
\rightarrow \gk,
\end{equation}
where $e_{\sigma}$ satisfies the ``functional equation'':
\begin{equation}
\label{fec}
e_{\sigma}(gp)=\chi (p)^{-1}e_{\sigma}(g), g\in(G)^U, p\in K.
\end{equation}

For homogeneous holomorphic line bundles \cite{tirao} ($\gk =\gc$) the
 functions in eq. (\ref{fec}) are holomorphic.

2). {\it Borel-Weil theorem\/} (cf. Ref. \cite{sbw})
 For every irreducible
representation $\pi_j$ of dominant weight $j$ of the compact connected
semisimple Lie group $G$ corresponds on every homogenous K\"ahlerian space
$G/K\approx G^{\gcm}/P_j$   a complete linear system $|D|$.
The representation space $\gh_j$ of the representation $\pi_j$ is the dual of
$\gl (D)$. The associated line
bundle ${\bf M}'$ is ample iff  the  space $G/K\approx G^{\gcm}/P_j$
is strictly
associated to the representation $\tilde{\pi}_j$ . Here $G^{\gcm} $ is the
complexification of $G$, $P_j$ is the parabolic group corresponding to the
dominant weight $j$ of the representation $\pi_j$ and
$|D|\approx{\bf P}(\gl (D))$.

3). The following theorem summarises some properties of {\it flag manifolds}
with significance for the present paper \cite{wolf}.
Let $X=G^{\gcm} /P$ be a complex manifold, where
$G^{\gcm}$ is a complex semisimple Lie group and $P$
is a parabolic subgroup. The following conditions are equivalent:
a)  $X=G^{\gcm}/P$ is compact;
b) $X$  is a compact simply connected \Ka~ manifold;
c) $X$ is a projective variety;
d) $X$ is a closed $G^{\gcm}$-orbit in a projective representation;
e) $X$ is a Hodge manifold and all homogeneous Hodge manifolds
are of this type.

4). A holomorphic line bundle ${\bf M}'$ on a compact complex manifold  \men~
is said {\it very ample} \cite{ss} if:
a) the set of divisors is without base points, i.e. there exists
a finite set of global sections $s_1,\ldots ,s_N\in \Gamma (\men, \bf{M'})$
such that for each $m\in \men$ at least one $s_j(m)$ is not zero;
b) the holomorphic map $\iota_{{\bf M}'}:\men\hookrightarrow\gcp^{N-1}$
given by
\begin{equation}
\iota_{{\bf M}'}=[s_1(m),\ldots ,s_N(m)]\label{scufund}
\end{equation}
is a holomorphic embedding.

The line bundle ${\bf M}'$ is said to be {\it ample} if there exists a positive
integer $r_0$ such that ${\bf M}'^r$ is very ample for all $r\geq r_0$. Note
that if ${\bf M'}$ is an ample line bundle on \men , then \men ~must be
projective-algebraic by Chow's theorem, hence \men~ is \Ka .

The concepts of ampleness and positivity for line bundles coincide. The
following theorem summarises the properties of ample line bundles
that are needed in this paper \cite{ss,bott,gh}. Below $[1]$ denotes
the hyperplane
line bundle on  the projective space \gph~  and the $C$-spaces are the simply
connected compact homogeneous manifolds.

 Let ${\bf M}'$ be a holomorphic line bundle  on a compact
complex manifold \men~of complex dimension $D$. The following conditions
are equivalent:
a) ${\bf M}'$ is positive;
b) the zero section of ${\bf M}'^*$ can be blown down to a point;
c) for all coherent analytic
sheaves ${\cal S}$ on \men~ there exists a
positive integer $m_0(\cal{S})$ such that
$H^i(\men,{\cal S}\otimes {\bf M}'^m)=0$ for $i>0,~m\geq m_0(\cal{S})$
(the vanishing theorem of Kodaira);
d) there exists a positive integer $m_0$ such that for all $m\geq m_0$,
there is an embedding  $\iota_{\bf M}:\men \hookrightarrow\gcp^{N-1}$ for
some $N\geq D$ such that ${\bf M}={\bf M'}^m$ is projectively induced, i.e.
${\bf M}=\iota ^*[1]$;
e) \men ~ is a Hodge manifold  (the embedding theorem of Kodaira);
f) the fundamental two-form of \men , the curvature matrix  and the first
Chern class of ${\bf M}'$ are related by the relations
$\omega=\sqrt{-1}/2\Theta_{{\bf M}'},
~c_1({\bf M}')=\omega/\pi$;
g) moreover, if \men~ is a K\"ahlerian $C$-space, then  \men~ is a flag
manifold.

5). We shall be concerned with  manifolds \men~ which admit an
embedding in some projective Hilbert space
\begin{equation}\label{emb}
\iota :
\widetilde {\bf M} \hookrightarrow \gph .
\end{equation}

In this paper we shall restrict ourselves to {\it biholomorphic embeddings}
  $\iota$.
Because $\iota$ in formula (\ref{emb}) is injective and holomorphic,
then it is a {\it k\"ahlerian embedding}, i.e.
\begin{equation}\label{pback}
\omega_{\men}=\iota^*\omega_{\gph},
\end{equation}
where $\omega$ is the fundamental two-form (i.e. closed,
(strongly)
non-degenerate)  of the K\"ahler manifold   and $\iota^*$ is the pull-back
of the mapping $\iota$.  Equivalently,  $\iota$ is an isometric embedding.

\section{\hspace{-4mm}.\hspace{2mm}THE COHERENT STATE AND COHERENT VECTOR
MANIFOLDS}

Let $ \xi :\gh^{\star}=\gh \setminus \{0\}\rightarrow\gph,
~ \xi (\mb{z})=[\mb{z}]$
be the mapping  which associates to the point $\mb{z}$ in the punctured
Hilbert space the linear
subspace $[\mb{z}]$ generated by $\mb{z}$, where $[\lambda\mb{z}]=[\mb{z}],
\lambda\in\gc^{*}$. The hermitian scalar product $(\cdot,\cdot )$ on \gh~
 is linear in
the second argument.

 Let us consider the principal bundle (\ref{pbundle}) and let us suppose
the existence of a map $e:G\rightarrow \gh^*$  as in eq. (\ref{sectiuni})
with the property (\ref{fec}) but globally
defined, i.e. on the neighbourhood (\ref{local}) $(G)^{\men}$.
Then $e(G)$ is called {\it  family of coherent vectors\/} \cite{raw}.
If there is a morphism of principal bundles \cite{hus},
i.e. the following diagram is commutative,

\begin{eqnarray}
\label{coherent}
\nonumber G &
~\stackrel{e}{\longrightarrow}~
& \gh^\star \\
\lambda\downarrow &  & \downarrow\xi\\ \label{comdiag}
\nonumber\men & ~\stackrel{\iota}{\longrightarrow}~ & \gph
\end{eqnarray}
then $\iota (\men )$  is called {\it family of coherent states corresponding
to the family of coherent vectors} $e(G)$ \cite{raw}.

We restrict ourselves to the case where the mapping $\iota$  is an embedding
of the homogeneous manifold \men~\cite{per}. The manifold \men ~ is called
{\it coherent state manifold} and the
$G$-homogeneous line bundle ${\bf M}_{\chi}$ is called {\it coherent vector
manifold} \cite{sbcag}.

Let now $\tilde{\pi}$ be a projective (in physical literature \cite{bargmann}
``ray'') representation  associated to the unitary irreducible representation
$\pi$ and $\tilde{G}$ the group of transformations which leaves invariant the
transition probabilities in the complex separable Hilbert space $\gh$.
If we use the projection
$\xi '=\xi_{|\sfera}$, i.e. $\xi ' :\sfera\rightarrow\gph,~
\xi '(\psi )=\tilde{\psi}=\{e^{i\phi}\psi|\phi\in\gr\}$, where \sfera~ is the
unit sphere in \gh , then $\tilde{\pi}\circ\xi '=\xi'\circ \pi$. The
triplet $(\tilde{\pi}, \tilde{G}, \gh )$ is a quantum system with symmetry in
the sense of Wigner and Bargmann \cite{wigner,bargmann}.
Then the manifold $\men\approx G/K$
can be realized as the orbit $\men =\{\tilde{\pi}(g)\tilde{\mb{e}}_0|g\in G\}$,
where $K$ is the stationary group of $\tilde{\mb{e}}_0$ and
$\mb{e}_0\in \gh^{\star}$ is fixed.
For a compact connected simply connected Lie group $G$,  the
existence of the representation $\tilde{\pi}$ implies the existence of the
unitary irreducible representation $\pi$ (cf. the theorem of
Wigner and Bargmann \cite{wigner,bargmann}). This implies the
existence of cross sections $\sigma :\men\rightarrow \sfera$. However, the
(Hopf) principal bundle $\mb{\xi}'=(\sfera ,\xi ', \gph )$ is a $U(1)$-bundle
and
in the construction of coherent vector manifold we need line bundles. But the
principal line bundle $\mb{\xi}'$ is obtained from the (tautological)
line bundle $[-1]=\mb{\xi}=(\gh^\star , \xi ,\gph )$ reducing the group
structure from $\gc^*=GL(1,\gc )$ to $U(1)$.

Here we also stress that the theorem of Wigner and Bargmann  is essentially
\cite{emch} the (first) fundamental theorem of projective geometry
\cite{artin}.

In order to have the physical interpretation of the ``classical system''
obtained by dequantizing the quantum one \cite{cval,sbaa},
we restrict to \Ka~ manifolds \men .
For example, for a compact connected simply connected Lie group $G$,
$\men\approx G/K\approx G^{\gcm}/P$ is a  \Ka~ manifold and the Borel-Weil
theorem assures the  geometrical realisation of the representation $\pi_j$
and of the representation  space $\gh_j$ if $\mb{e}_0=j$.

The representation $\pi_j$ can be uniquely extended to the group homomorphism
$\pi^*_j:~ G^{\gcm}\rightarrow\pi_j^*(G^{\gcm})$, and respectively, Lie algebra
isomorphism
$\stackrel{.}{\pi}^*_j:\got{g}^{\gcm}\rightarrow\stackrel{.}{\pi}^*_j
(\got{g}^{\gcm})$ by
\begin{equation}
\pi^*_j(\exp (Z))=\mbox{e}^{\stackrel{.}{\pi}^*_j(Z)}, Z\in\got{g}^{\gcm},
\end{equation}
where ${\exp}:\got{g}^{\gcm}\rightarrow G^{\gcm}$ and $\mbox{e}
:\stackrel{.}{\pi}^*_j\rightarrow\pi^*_j(G)$ are exponential maps, while
$\stackrel{.}{\pi}^*_j$ is the complexification of the Lie algebra
$\stackrel{.}{\pi}_j(\got{g})$. We use also the notation $F_{\alpha}=
\stackrel{.}{\pi}^*_j(f_{\alpha})$,
where $\alpha $ is in the set $\Delta$ of the
roots of the Lie algebra $\got{g}$ of $G$ with generators $f_{\alpha}$
of the Cartan-Weyl base of $\got{g}^{\gcm}$ (see also \cite{sbcag}).

Then $\mb{e}_g \df e(g)\df \pi^* (g)\mb{e}_0, g\in G^{\gcm}$ is the
family of coherent vectors, while $\{\tilde{\mb{e}}\}_{g\in G^{\gcm}}$
is  the family of coherent states. The relation
$\mb{e}_g={\mbox{e}}^{i\alpha(g)}\mb{e}_{\lambda (g)}$
defines a  fibre bundle with base \men~ and fibre $U(1)$
\cite{per}. More precisely,
{\it the function}
\begin{equation}
\label{functiune}
\mb{\Upsilon} (g)=(\Upsilon,\mb{e}_g )
\end{equation}
{\it      is holomorphic on $G^{\gcm}$ and defines
holomorphic sections
on the homogeneous holomorphic  line bundle ${\bf M}'$ associated
to the principal line bundle $P\rightarrow G^{\gcm}\rightarrow G^{\gcm}/P$ by
the holomorphic character} $\chi$
\begin{equation}
\label{car}
\pi^* (p)e_0=\chi^{-1}(p)e_0,~p\in P,
~\chi (p)={\mbox{\rm e}}^{-i\alpha (p)}.
\end{equation}
Indeed, the function $\mb{\Upsilon} (g)$ verifies
$\mb{\Upsilon} (gp)=\chi^{-1}(p)\mb{\Upsilon} (g), g\in G^{\gcm},~p\in P$,
i.e. eq. (\ref{fec}),
and the corresponding  holomorphic sections are associated via
eq. (\ref{sectiuni}).

Let also  the function
\begin{equation}
\label{u1}
\mb{\Upsilon} '(g)=\mb{\Upsilon} '(gP):=
\displaystyle{\frac{\mb{\Upsilon} (g)}
{(\mb{e}_0 ,\mb{e}_g)}},
\end{equation}
defined on the set
\begin{equation}\label{nonzero}
(\mb{e}_0 ,\mb{e}_g)\not= 0.
\end{equation}
Then
\begin{equation}
\label{uu2}
\mb{\Upsilon}':{\cal V}_0\rightarrow\gc,
~\mb{\Upsilon} '(Z)=(\Upsilon ,\mb{e}_{Z,j}),
\end{equation}
where the  Perelomov's coherent vectors  are
\begin{equation}
{\mb e}_{Z,j}=\exp\sum_{{\varphi}\in\Delta^+_n}(Z_{\varphi}F^+_{\varphi})
j ,~~~~\underline{{\mb e}}_{Z,j}=({\mb e}_{Z,j} ,
{\mb e}_{Z,j})^{-1/2}{\mb e}_{Z,j} ,
\label{z}
\end{equation}
\begin{equation}
{\mb e}_{B,j}=\exp\sum_{{\varphi}\in\Delta^+_n}(B_{\varphi}F^+_{\varphi}-{\bar
B}_{\varphi}F^-_{\varphi}) j , \label{b}
~~~{\mb e}_{B,j}\df  \underline{{\mb e}}_{Z,j}.
\end{equation}
Here $\Delta^+_n$ denotes the positive non-compact roots, $ Z\df (Z_ \varphi )
\in {\gc}^D$  are local
coordinates in the maximal neighbourhood
${\cal V}_0 \subset \widetilde {\bf M} $.
In eqs. (\ref{z}), (\ref{b}) $F^+_{\varphi} j\neq  0,F^-_{\varphi} j = 0,
~  \varphi\in\Delta^+_n $.

The system \{$\mb{e}_g$\}, $g\in G^{\gcm}$ is overcomplete
\cite{berezin,per,onofri}
and $(\mb{e}_g,\mb{e}_{g'})$, up to a
factor, is a reproducing kernel for the holomorphic vector
bundle $\xi_0:{\bf M}\rightarrow\men$ \cite{zb}.

\section{\hspace{-4mm}.\hspace{2mm} CUT LOCUS AND COHERENT STATES}

Let $X$ be complete Riemannian manifold. The point  $q$ is in the
{\it  cut locus} ${\bf CL_p}$ of $p\in X$ if $q$ is the nearest point to $p$ on
the geodesic emanating from $p$ beyond which the geodesic ceases to minimize
his arc length (cf. \cite{kn}, see also Ref. \cite{ber2} for more references).

\begin{rem} $ {\mbox{\rm codim}}_{\gcm} {\bf CL}_p \geq 1$.
\label{rem4}
\end{rem}

We call {\it  polar divisor} of $\mb{e}_0$ the set $
\Sigma _0=\left\{ \mb{e} \in e(G)|(\mb{e}_0, \mb{e})=0\right\}$. This
denomination is inspired after Wu \cite{wu}, who used this term in the case
of the complex Grassmann manifold \Gras .

Let ${\got g}={\got k}\oplus {\got m}$ be  the orthogonal
decomposition of {\got g} with respect to the $B$-form,
${\rm Exp}_p:\widetilde {\bf M}_p\rightarrow \widetilde {\bf M}$ the
the geodesic exponential map and
 $o=\lambda (e)$, where $e $ is the unit element
in $G$. Then $\got{m}$ is identified with the tangent space at $o$, $\men _o$,
and $\men\approx \exp\got{m}$.

Let us consider the following two conditions

\parbox{1cm}{\it A)}\parbox[t]{132mm}{$ {\rm Exp} _o
=\lambda \circ \exp \vert _{\got m}~.  $}

\parbox{1cm}{\it B)}\parbox[t]{132mm}{ On the Lie algebra  {\got g} of $ G$
there exists an $ Ad(G)$-invariant,
symmetric, non-degenerate bilinear form $ B$ such that the restriction
of $B$ to the Lie algebra {\got k}  of $ K$  is likewise non-degenerate.}

Note that {\it the symmetric spaces have property $A)$} and if $\men \approx
G/K$ {\it verifies} $B)$, {\it then it also verifies} $A)$ (cf. \cite{kn}).

In \cite{ber2} it was proved the following

\begin{thm}
  Let \men~    be a homogeneous manifold
$\widetilde {\bf M}\approx G/K$.  Suppose that there exists a unitary
irreducible representation $\pi_j$ of $G$ such that in a neighbourhood
$ {\cal V}_0$  around $Z=0$ the coherent states are parametrized as
in eq. (\ref{z}).  Then the manifold \men~ can be represented as
the disjoint union

\begin{equation}
\widetilde {\bf M} ={\cal V}_0\cup \Sigma _0.\label{reu}
\end{equation}

Moreover, if the condition  $B)$ is true, then
\begin{equation}
\Sigma _0={\bf CL}_0. \label{clo}
\end{equation}
\end{thm}

\begin{cor}
 Suppose that $\widetilde {\bf M}$ verifies  $ B)$  and admits the embedding
(\ref{emb}). Let $0, Z\in
\men $.  Then  $Z\in {\bf CL}_0$ iff  the Cayley distance
between the images  $\iota (0), \iota (Z)\in \gph$ is  $\pi /2$
\begin{equation}
d_c(\iota (0),\iota (Z))=\pi /2.
\end{equation}
\end{cor}

Here $d_c$ denotes the the hermitian elliptic Cayley distance on the projective
space
\begin{equation}
d_c([\omega '],[\omega ])=\arccos \frac{\vert(\omega ',\omega)\vert}
{\Vert \omega '\Vert \Vert \omega \Vert }~. 
\end{equation}

\section{\hspace{-4mm}.\hspace{2mm} DISCUSSION}

1). We now state precisely the consequences of theorem 1
for Berezin's quantization on symmetric compact complex manifolds
(see firstly the second reference in \cite{berezin}).

 Expressing the  supercompleteness  of the system of coherent vectors
$\{\mb{e}_g\}$
by the Parseval identity, Berezin  introduces the Hilbert
space ${\cal F}_h$ of holomorphic functions  (\ref{uu2}) on
${\cal V}_0$, denoted $\mb{\Upsilon}(z)$, with the scalar product
\begin{equation}
(\mb{\Upsilon}_1,\mb{\Upsilon}_2)
=\tilde{c}(h)\int_{{\cal V}_0}\mb{\Upsilon}_1(z)\overline{\mb{\Upsilon}}_2(z)
F^{\frac{1}{h}}(z,\overline{z})d\mu (z,\overline{z}),\label{berezin}
\end{equation}
Here
$-\ln F$ is the K\"ahler potential on \men ,
$d\mu (z,\overline{z})=\pi ^{-n}F(z,\overline{z})
d\mu_L(z,\overline{z})$, and $\mu_L$ denotes the Lebesgue measure on \men .
In Berezin's terminology,
$h$ belongs to the {\it admissible set}, $z$ are {\it special coordinates}
on ${\cal V}_0\subset\men$ and $z=0$  is a {\it distinguished point}.

Here we just comment that, due to Remark \ref{rem4}, {\it the integration in}
eq. (\ref{berezin}) {\it can be extended to all} \men  .

In Berezin's formulation,
{\it the quantization algebra \got{A}}, which is a special quantization
with correspondence principle in the weak form, {\it is restricted also only
on} ${\cal V}_0$. We remember Berezin's definition of the $*$-product of
covariant symbols
\begin{equation}\label{star}
(A_1*A_2)(z,\overline{z})=
\int A_1(z,\overline{v})A_2(v,\overline{z})
G_h(z,\overline{z}|v,\overline{v})d\mu (v,\overline{v}),
\end{equation}
where the kernel and the covariant symbol attached to the
operator $\hat{A}$ are, respectively
\begin{equation}
\label{ker1}
G_h(z,\overline{z}|v,\overline{v})=c(h)\displaystyle{
\frac{\mb{\Upsilon}_{\overline{v}}(z)\mb{\Upsilon}_v(\overline{z})}
{\mb{\Upsilon}_{\overline{z}}(z)\mb{\Upsilon}_{\overline{v}}(v)}},
\end{equation}
\begin{equation}
\label{a1}
A(z,\overline{v})=\displaystyle{\frac{(\hat{A}\mb{\Upsilon}_{\overline{v}},
\mb{\Upsilon}_{\overline{z}})}
{(\mb{\Upsilon}_{\overline{v}},\mb{\Upsilon}_{\overline{z}})}}.
\end{equation}
Here $\mb{\Upsilon}_{\overline{v}}(z)=L_h(z,\overline{v})=
F^{-1/h}(z,\overline{v})$ and $(f,\mb{\Upsilon}_{\overline{v}})=f(v)$.
 If $f_k(z)$ is an orthonormal basis in ${\cal F}_h$, then
$L_h(z,\overline{v})=\sum f_k(z)\overline{f}_k(v),$
and $L_h$ is the kernel of the (Bergman) projector
$P_B:L^2(F^{1/h}d\mu)\rightarrow {\cal F}_h: (P_B)f(z)=\int
L_h(z,\overline{\xi})f(\xi ,\overline{\xi})d\mu(\xi,\overline{\xi})$.
{\it We can  explicitly write down  the
domain of definition of the functions} $\mb{\Upsilon}_{\overline{v}}(z)
\in{\cal F}_h$,
: {\it at fixed} $v,~z\notin {\bf CL_v}.$

2). Referring to  Onofri's paper \cite{onofri} (see also the papers of
Rawnsley \cite{raw} and Moreno \cite{moreno}, we also comment that,
due to eqs. (\ref{nonzero}), (\ref{uu2}) and (\ref{reu}),
{\it the functions expressing the
coherent states are not defined on all the manifold} and their domain of
definition has the  geometrical significance given by the relation
(\ref{clo}) in theorem 1.

3).  Cahen, Gutt and Rawnsley \cite{cgr} formulated globally Berezin's
construction of covariant symbols of operators in terms of sections of
the prequantization
bundle $({\bf M},h,\nabla )$  of the \Ka~ manifold $(\men,\omega)$.
Berezin's definition (\ref{star}) of the
$*$-product is modified by the presence of the function $\epsilon$  as
\begin{equation}\label{star1}
(A_1*_kA_2)(x)=\int_{\men}A_1(x,y)A_2(y,x)\Psi^k(x,y)
\epsilon ^{(k)}\displaystyle{\frac{k^n\omega^n}{n!}}.
\end{equation}

We remember the notation.
Let $q\in {\bf M}_x=\xi^{-1}_0(x)$ be a fixed frame field over
\nolinebreak\men~
and the holomorphic section $s\in\Gamma (\men ,{\bf M})$.
The evaluation of section at $x$ gives $s(x)=\mb{\Upsilon}_q(s)q$
and  the continuous coordinate function corresponds to
$\mb{\Upsilon}\in{\cal F}_h$ in Berezin's notation.
The unique element $e_q\in\gh$ determined by  Riesz theorem from the relation
$\mb{\Upsilon}_q(s)=(s,e_q)$ verifies the definition (\ref{coherent}),
(\ref{fec})
of coherent states with \nolinebreak$\chi (c)=c$.\nopagebreak

The Berezin's symbol in eq. (\ref{star1}) is defined in terms of
sections (its analogue for functions is given by eq. (\ref{a1}))
under the restriction (\ref{nonzero}),  ($(e_{q'},e_q)\not= 0)$:
\begin{equation}
\label{symbol2}
A(x,y)=\displaystyle{\frac{(\hat{A}e_{q'},e_q)}{(e_{q'},e_q)}},~
\xi_0(q)=x,~\xi_0(q')=y.
\end{equation}
So, {\it the Berezin symbol} $A(x,y)$  in eq. (\ref{symbol2})
{\it is defined at a fixed} $x$ {\it for} $y\notin {\bf CL}_x$.

{\it The characteristic function} is  defined in Ref. \cite{cgr}
in a neighbourhood $U\times\overline{U}$ of the diagonal of
$\men\times\overline{\men}$ as
\begin{equation}
\label{pcs}
\tilde{\Psi} (x,y)=\displaystyle{\frac{|s|^2(x,\overline{x})
|s|^2(y,\overline{y})}{||s|^2(x,\overline{y})|^2}},
\end{equation}
and it is related to the Calabi's diastasis \cite{calabi}  by the relation
$D=-2\log \tilde{\Psi}.$
In eq. (\ref{pcs}) $|s|^2(x,\overline{x})=h_x(s(x),s(x))$ and
$|s|^2(x,\overline{y})$ is its analytic continuation in
$U\times\overline{U}$. In fact, due to theorem 1,
{\it $\tilde{\Psi}(x,y)$ and $D(x,y)$ are defined at a fixed $x$ for} $y\notin
{\bf CL}_x$.

The 2-{\it point function}, whose local analogue is given by eq. (\ref{ker1}),
is globally defined in terms of sections as
\begin{equation}
\label{psi2}
\Psi (x,y)=\displaystyle{\frac{|(e_{q'},e_q)|^2}
{||e_{q'}||^2||e_q||^2}},~\xi_0 (q)=x,~
\xi_0(q')=y .
\end{equation}
If the quantization is regular,     i. e. $\epsilon =$ ct,
then $\tilde{\Psi}= \Psi$ and eq. (\ref{star1}) can be put in the form
(\ref{star}).  We remember also that $\epsilon (x)=||e_q||^2h(q,q),
\xi_0 (q)=x, x\in\men$.

4). Finally, let me mention that during the Workshop  I have received
a positive answer \cite{martin}
to my question: {\it is the polar divisor in theorem 1 a divisor
in the sense of algebraic geometry?\/} A more precise statement
relative to the advanced question can be given. Indeed, let $[~]$ be the
functorial
homomorphism  $[~]: \mbox{Div}(\men )\rightarrow H^1(\men ,{\cal O}^*)$
between the group of divisors of a complex manifold and the
Picard group of equivalence classes of $C^{\infty}$ line bundles \cite{gh}.
Then \cite{project} :

{\it Let \men~ be a simply connected Hodge manifold admitting the embedding
(\ref{emb}). Let ${\bf M}=\iota^{*}[1]$ be the
unique, up to equivalence, projectively induced line bundle with a given
admissible connection. Then ${\bf M}=[\Sigma_0]$. Moreover, if the homogeneous
manifold \men~ verifies condition B), then ${\bf M} =[{\bf CL}_0]$.
In particular,
the first relation is true for \Ka ian $C$-spaces, while the second one for
hermitian symmetric spaces\/.}

{\bf Acknowledgments}

The author expresses his thanks to Professors
Anatol Odzijewicz and Aleksander  Strasburger for the possibility
to attend  the XV-th  Workshop on Geometric Methods in Physics.
Discussions during the Workshop with Professors Michel Cahen, Gerard Emch,
Simon Gutt,
Joachim Hilgert, and Martin Schlichenmaier are kindly acknowledged.

\end{document}